\documentclass[11pt]{article}
\usepackage{amsmath}
\usepackage{amssymb}
\usepackage{amsthm}

\newcommand\nonu{\nonumber}
\newcommand\ma{\medskipamount}
\newcommand\ba{\bigskipamount}
\newcommand\sa{\smallskipamount}
\newcommand\sLP{\\[\sa]}
\newcommand\bLP{\\[\ba]}
\newcommand\mPP{\\[\ma]\indent}
\newcommand\bPP{\\[\ba]\indent}

\newcommand\CC{\mathbb{C}}
\newcommand\RR{\mathbb{R}}
\newcommand\FSA{{\cal A}}
\newcommand\FSB{{\cal B}}
\newcommand\FSH{{\cal H}}
\newcommand\al\alpha
\newcommand\be\beta
\newcommand\de\delta
\newcommand\tha\theta
\newcommand\la\lambda
\newcommand\om\omega
\newcommand\Ga{\Gamma}

\newcommand\half{\frac12}
\newcommand\thalf{\tfrac12}
\newcommand\iy\infty
\newcommand\id{\operatorname{id}}
\newcommand\lan{\langle}
\newcommand\ran{\rangle}
\newcommand\Equiva{\Longleftrightarrow}
\newcommand{\hyp}[5]{\,\mbox{}_{#1}F_{#2}\!\left(
  \genfrac{}{}{0pt}{}{#3}{#4};#5\right)}
\newcommand\LHS{left-hand side}
\newcommand\RHS{right-hand side}
\newcommand\Span{\operatorname{Span}}

\newcounter{alphc} 
\newcommand\alphlist{\begin{list}{\alph{alphc})}{\usecounter{alphc} 
 \parsep0cm \itemsep0.2cm \topsep0.2cm}} 
\newcommand\bma{\begin{pmatrix}}
\newcommand\ema{\end{pmatrix}}

\numberwithin{equation}{section}
\newtheorem{theorem}{Theorem}[section]
\newtheorem{proposition}[theorem]{Proposition}
\newtheorem{lemma}[theorem]{Lemma}
\newtheorem{Definition}[theorem]{Definition}
\newenvironment{definition}{\begin{Definition}\rm}{\end{Definition}}
\newtheorem{Remark}[theorem]{Remark}
\newenvironment{remark}{\begin{Remark}\rm}{\end{Remark}}
\newtheorem{Example}[theorem]{Example}
\newenvironment{example}{\begin{Example}\rm}{\end{Example}}
\newcommand\Proof{\smallskip\noindent{\bf Proof}\quad}

\begin{document}
\title{Representations of SU(2) and Jacobi polynomials
\footnote{Slight revision and extension
of notes written for the last topic of
an honours class on {\em Orthogonal polynomials},
University of Amsterdam, January 2007}}
\author{Tom H. Koornwinder}
\date{}
\maketitle
\begin{abstract}
This is a tutorial introduction to the representation theory of SU(2)
with emphasis on the occurrence of Jacobi polynomials in the matrix
elements of the irreducible representations. The last section traces
the history of the insight that Jacobi polynomials occur in the
representation theory of SU(2).
\end{abstract}
%
%
\section{Introduction}
These are lecture notes, dating back to 2007, which present elements of the
representation theory of SU(2) and the occurrence there of Jacobi polynomials
to an audience of advanced undergraduate students in math. Some statements
are left as exercises to the reader. Some longer exercises are collected in Section~\ref{53}.
Some historical background is given in
in Section~\ref{52}.

Our treatment of the theory does not claim originality.
See generalities about the representation theory of SU(2) in
Sugiura \cite{3}. Some parts of books dealing with
representation theory of SU(2) including special function aspects
are Vilenkin \cite[Ch.~III]{12},
Miller \cite[Section 5.16]{13},
Biedenharn \& Louck \cite[Section 3.6]{8},
Vilenkin \& Klimyk \cite[Ch.~6]{4}, and
Andrews, Askey \& Roy \cite[Sections 9.11--9.16]{2}.
\section{Preliminaries about representation theory}
Let $G$ be a group. Representations of $G$ can be defined on any vector space
(possibly infinite dimensional) over any field, but we will only
consider representations on finite dimensional complex vector spaces.
Let $V$ be a finite dimensional complex vector space.
Let ${\rm GL}(V)$ be the set of all invertible linear transformations of $V$.
This is a group under composition. If $V$ has dimension $n$ and if we choose
a basis $e_1,\ldots,e_n$ of $V$ then the map
$x=x_1e_1+\cdots+x_ne_n\mapsto(x_1,\ldots,x_n)\colon V\to\CC_n$ is an
isomorphism of vector spaces. There is a corresponding group isomorphism
${\rm GL}(V)\to {\rm GL}(\CC^n)$ which sends each invertible
linear transformation
of $V$ to the corresponding invertible matrix with respect to this basis.
We denote ${\rm GL}(\CC^n)$ by ${\rm GL}(n,\CC)$: the group of all invertible
complex $n\times n$ matrices. Here the group multiplication is by
multiplication of matrices.
\begin{definition}
A {\em representation}
of a group $G$ on a finite dimensional complex vector space
$V$ is a group homomorphism $\pi\colon G\to {\rm GL}(V)$.
A linear subspace $W$
of $V$ is called {\em invariant} (with respect to the representation $\pi$)
if $\pi(g)\,W\subset W$ for all $g\in G$.
The representation $\pi$ on $V$ is called {\em irreducible} if $V$ and
$\{0\}$ are the only invariant subspaces of $V$.
\end{definition}
\begin{definition}
Let $\pi$ be a representation of a group $G$ on a finite dimensional
complex vector space $V$. Choose a basis $e_1,\ldots,e_n$ of $V$.
Then, for $g\in G$, the linear map $\pi(g)$ has a matrix
$(\pi_{i,j}(g))_{i,j=1,\ldots,n}$ with respect to this basis, which is
determined by the formula
\begin{equation*}
\pi(g)\,e_j=\sum_{i=1}^n \pi_{i,j}(g)\,e_i.
\end{equation*}
The $\pi_{i,j}$ are complex-valued functions on $G$ which are called the
{\em matrix elements} of the representation $\pi$ with respect to the basis
$e_1,\ldots,e_n$.
\end{definition}
\begin{remark}
Let ${\rm End}(V)$ be the space of all linear transformations $A\colon V\to V$.
If $\pi$ is a map of the group $G$ into ${\rm End}(V)$ such that
$\pi(g_1g_2)=\pi(g_1)\pi(g_2)$ for all $g_1,g_2$ and $\pi(e)=\id$, then
$\pi$ maps into ${\rm GL}(V)$ and $\pi$ is a representation of $G$ on $V$
(proof left as an exercise).
\end{remark}
\begin{definition}
A {\em topological group} is a set $G$ which is both a group and a
topological space such that the maps
$(g_1,g_2)\mapsto g_1g_2\colon G\times G\to G$ and
$g\mapsto g^{-1}\colon G\to G$ are continuous.
\end{definition}
\begin{example}
${\rm GL}(n,\CC)$ can be considered as a subset of $\CC^{n^2}$ by associating
with the element
$T=(t_{i,j})_{i,j=1,\ldots,n}\in {\rm GL}(n,\CC)$ the $n^2$
complex coordinates
$t_{i,j}$. Then the group ${\rm GL}(n,\CC)$, with the topology inherited
from
$\CC^{n^2}$, is a topological group
(proofs as an exercise).

Let $V$ be an $n$-dimensional complex vector space.
With respect to any basis of $V$ the group ${\rm GL}(V)$ is isomorphic with
${\rm GL}(n,\CC)$.
Give a topology to ${\rm GL}(V)$ such that this isomorphism is also a homeomorphism.
Then ${\rm GL}(V)$ is a topological group and the topology is independent of
the choice of the basis
(proofs as an exercise).
\end{example}
\begin{definition}
A {\em representation}
of a topological group $G$ on a finite dimensional complex vector space
$V$ is a continuous group homomorphism $\pi\colon G\to {\rm GL}(V)$. 
\end{definition}
\begin{remark}
\label{3}
Let $G$ be a topological group, $V$ a finite dimensional complex vector space
and $\pi\colon G\to {\rm GL}(V)$ a group homomorphism.
Let $e_1,\ldots,e_n$ a basis for $V$.
Then the following five properties are equivalent:
\alphlist
\item
$\pi$ is continuous;
\item
for all $v\in V$ the map
$g\mapsto \pi(g)\,v\colon G\to V$ is continuous;
\item
for all $j$ the map $g\mapsto \pi(g)\,e_j\colon G\to V$ is continuous;
\item
for all  $v\in V$ and for all complex linear functionals $f$ on $V$\\
the map $g\mapsto f(\pi(g) v)\colon G\to\CC$ is continuous.
\item
The matrix elements $\pi_{i,j}$ of $\pi$ with respect to the basis
$e_1,\ldots, e_n$ are continuous functions on $G$.
\end{list}
The proofs are left as an exercise.
Be aware that these equivalences are not necessarily true if $V$ is an
infinite dimensional topological vector space.
\end{remark}
\begin{remark}
If $\pi$ is a representation of $G$ on $V$ and if $H$ is a subgroup of $G$
then the restriction of the group homomorphism $\pi\colon G\to {\rm GL}(V)$
to
$H$ is a group homomorphism $\pi\colon H\to {\rm GL}(V)$, so it is
a representation of $H$ on $V$.

If $G$ is moreover a topological group then $H$
with the topology inherited from $G$ becomes a topological group (exercise!).

If, furthermore, $\pi$ is a representation of $G$ as a topological group
on $V$ then the restriction of $\pi$ to $H$ is
a representation of $H$ as a topological group on $V$.
\end{remark}
\begin{definition}
\label{8}
Let $V$ be a finite dimensional complex vector space with hermitian inner
product $\lan\,.\,,\,.\,\ran$. A representation $\pi$ of a group $G$ on $V$ is
called {\em unitary} is $\pi(g)$ is a unitary operator on $V$ for all
$g\in G$, i.e., if
\begin{equation*}
\lan\pi(g)\,v,\pi(g)\,w\ran=\lan v,w\ran\quad
\mbox{for all $v,w\in V$ and for all $g\in G$.}
\end{equation*}
\end{definition}
\begin{remark}
\label{10}
Let $V$ and $G$ be as in Definition \ref{8} and let $\pi$ be a representation
of $G$ on $V$. Let $e_1,\ldots,e_n$ be an orthonormal basis of $V$ and
let $\pi(g)$ have matrix $(\pi_{i,j}(g))$ with respect to this basis.
Then the representation $\pi$ is unitary iff the matrix $(\pi_{i,j}(g))$
is unitary for each $g\in G$.
One of the ways to characterize unitarity of the matrix $(\pi_{i,j}(g))$
is that
\begin{equation*}
\overline{\pi_{i,j}(g)}=\pi_{j,i}(g^{-1})\quad(i,j=1,\ldots,n).
\end{equation*}
The proof is left as an exercise.
\end{remark}
\begin{proposition}
{\bf (Complete reducibility of unitary representations)}\\
Let $V$ and $G$ be as in Definition \ref{8} and let $\pi$ be a
unitary representation of $G$ on $V$. Then:
\alphlist
\item
If $W$ is an invariant subspace of
$V$ then the orthoplement $W^\perp$ of $W$ is also an invariant subspace.
\item
$V$ can be written as an orthogonal direct sum of subspaces $V_i$ such that
the representation $\pi$, when restricted to $V_i$, is irreducible.
\end{list}
\end{proposition}
The proof is left as an exercise.
\section{A class of representations of SU(2)}
\label{47}
Fix $l\in\{0,\thalf,1,\ldots\}$.
Let $\FSH_l$ be the space of homogeneous polynomials of degree $2l$ in two
complex variables $z_1,z_2$. So the monomials $z_1^{l-n}z_2^{l+n}$
($n=-l,-l+1,\ldots,l$) form a basis of $\FSH_l$, and $\FSH_l$ has dimension
$2l+1$. For reasons which will become clear later, we will work with a
renormalized basis
\begin{equation}
\psi_n^l(z_1,z_2):=\binom{2l}{l-n}^\half z_1^{l-n} z_2^{l+n}\qquad
(n=-l,-l+1,\ldots,l).
\label{1}
\end{equation}
For $A\in {\rm GL}(2,\CC)$ and $f\in \FSH_l$ define the function $t^l(A)f$
on $\CC^2$
by
\begin{equation}
(t^l(A)f)(z):=f(A'z)\qquad(z=(z_1,z_2)\in\CC^2),
\label{9}
\end{equation}
where $A'$ is the transpose of the matrix $A$. So
\begin{equation*}
\left(t^l\bma a&b\\c&d\ema f\right)(z_1,z_2)=f(az_1+cz_2,bz_1+dz_2),\quad
{\rm where}\; \bma a&b\\c&d\ema\in {\rm GL}(2,\CC).
\end{equation*}
From this it is clear that $(t^l(A)f)(z_1,z_2)$ is again a homogeneous
polynomial of degree $2l$ in $z_1,z_2$. Moreover, $t^l$ is a representation
of ${\rm GL}(2,\CC)$ on $\FSH_l$, since
$t^l(I)f=f$ and
\begin{multline*}
(t^l(AB)f)(z)=f((AB)'z)=f(B'A'z)=(t^l(B)f)(A'z)\\
=\bigl(t^l(A)(t^l(B)f)\bigr)(z)=\Bigl(\bigl(t^l(A)t^l(B)\bigl)f\Bigr)(z).
\end{multline*}

The matrix elements $t_{m,n}^l$ ($m,n=-l,-l+1,\ldots,l$) of $t^l$
with respect to the basis \eqref{1} are determined by
\begin{equation}
t^l(g)\,\psi_n^l=\sum_{m=-l}^l t_{m,n}^l(g)\,\psi_m^l\qquad(g\in
{\rm GL}(2,\CC)),
\label{32}
\end{equation}
Since
\begin{equation}
\left(t^l\bma a&b\\c&d\ema \psi_n^l\right)(z_1,z_2)=
\binom{2l}{l-n}^\half (az_1+cz_2)^{l-n}\,(bz_1+dz_2)^{l+n},
\label{31}
\end{equation}
\eqref{32} can be written more explicitly as
\begin{multline}
\binom{2l}{l-n}^\half (az_1+cz_2)^{l-n}\,(bz_1+dz_2)^{l+n}=
\sum_{m=-l}^l \binom{2l}{l-m}^\half\,t_{m,n}^l\bma a&b\\c&d\ema\,
z_1^{l-m}z_2^{l+m},\\
\bma a&b\\c&d\ema\in {\rm GL}(2,\CC).
\label{2}
\end{multline}
From \eqref{2} we see that $t_{m,n}^l\bma a&b\\c&d\ema$ is a homogeneous
polynomial of degree $2l$ in $a,b,c,d$, so $t_{m,n}^l$ is continuous on
${\rm GL}(2,\CC)$. By Remark \ref{3} $t^l$ is then also a representation of
${\rm GL}(2,\CC)$ considered as a topological group.

For fixed $n$ we can consider \eqref{2} as a {\em generating function}
for the matrix elements $t_{m,n}^l$ with $m=-l,\ldots,l$: the matrix elements
are obtained as the coefficients in the power series expansion of the
elementary function in $z_1,z_2$ on the \LHS.

From \eqref{2} for $n=l$ elementary expressions for the matrix elements
$t_{n,l}^l$ can be obtained (exercise):
\begin{equation}
t_{m,l}^l\bma a&b\\c&d\ema=\binom{2l}{l-m}^\half b^{l-m} d^{l+m}.
\label{35}
\end{equation}

From \eqref{2} we can derive a double generating function for the
matrix elements $t_{m,n}^l$:
Multiply both sides of \eqref{2} with
\begin{equation*}
\binom{2l}{l-n}^\half w_1^{l-n} w_2^{l+n},
\end{equation*}
and sum over $n$. Then we obtain
\begin{multline}
(az_1w_1+bz_1w_2+cz_2w_1+dz_2w_2)^{2l}=
\sum_{m,n=-l}^l \binom{2l}{l-m}^\half \binom{2l}{l-n}^\half\,
t_{m,n}^l\bma a&b\\c&d\ema\\
\times z_1^{l-m}z_2^{l+m}\,w_1^{l-n}w_2^{l+n},\qquad
\bma a&b\\c&d\ema\in GL(2,\CC).
\label{4}
\end{multline}

Formula \eqref{4} implies the symmetry
\begin{equation}
t_{m,n}^l\bma a&b\\c&d\ema=t_{n,m}^l\bma a&c\\b&d\ema,
\label{5}
\end{equation}
while \eqref{2} implies that
\begin{equation}
t_{m,n}^l\bma a&b\\c&d\ema=t_{-m,-n}^l\bma d&c\\b&a\ema.
\label{6}
\end{equation}
From \eqref{5} and \eqref{6} we obtain a third symmetry
\begin{equation}
t_{m,n}^l\bma a&b\\c&d\ema=t_{-n,-m}^l\bma d&b\\c&a\ema.
\label{7}
\end{equation}
The details of the proofs of \eqref{5}--\eqref{7} are left as exercises.

Let SU(2) denote the set of all $2\times 2$ unitary matrices of
determinant 1. This is clearly a subgroup of ${\rm GL}(2,\CC)$ (exercise!).
Note that SU(2) consists of all matrices
\begin{equation}
\bma a&-\overline c\\c&\overline a\ema\quad
\mbox{with $a,c\in\CC$ and $|a|^2+|c|^2=1$.}
\label{20}
\end{equation}
Prove this as an exercise. Hence, as a topological space, SU(2) is
homeomorphic with $\{(a,c)\in\CC^2\mid |a|^2+|c|^2=1\}$, which is the
unit sphere in $\CC^2$, i.e., the sphere $S^3$. In particular,
SU(2) is compact.

The representation $t^l$ of ${\rm GL}(2,\CC)$ given by \eqref{9}, becomes by
restriction a representation of SU(2). Put a hermitian inner product on
$\FSH_l$ such that the basis of functions $\psi_n^l$
($n=-l,-l+1,\ldots,l$) is orthonormal.
\begin{proposition}
The representation $t^l$ of $\rm SU(2)$ is unitary.
\end{proposition}
\Proof
The inverse of
$\bma a&-\overline c\\c&\overline a\ema\in {\rm SU(2)}$ is
$\bma \overline a&\overline c\\-c&a\ema$.
In view of Remark \ref{10} we have to show that
\begin{equation*}
\overline{t_{m,n}^l\bma a&-\overline c\\c&\overline a\ema}
=t_{n,m}^l\bma \overline a&\overline c\\-c&a\ema.
\end{equation*}
Since, by \eqref{2}, $t_{m,n}^l\bma a&b\\c&d\ema$ is a polynomial with
real coefficients in $a,b,c,d$, we have
\begin{equation*}
\overline{t_{m,n}^l\bma a&-\overline c\\c&\overline a\ema}=
t_{m,n}^l\bma \overline a&-c\\\overline c&a\ema.
\end{equation*}
Hence we have to show that
\begin{equation*}
t_{m,n}^l\bma \overline a&-c\\\overline c&a\ema=
t_{n,m}^l\bma \overline a&\overline c\\-c&a\ema.
\end{equation*}
This last identity follows from \eqref{5}.
\qed
\section{Computation of matrix elements of representations of SU(2)}
We can use the generating function \eqref{2} in order to compute
the matrix elements $t_{m,n}^l$. First we expand the two powers on the
\LHS\ of \eqref{2} by the binomial formula:
\begin{align*}
(az_1+cz_2)^{l-n}&=
\sum_{j=0}^{l-n}\binom{l-n}j a^j z_1^j c^{l-n-j} z_2^{l-n-j},\\
(bz_1+dz_2)^{l+n}&=
\sum_{k=0}^{l+n}\binom{l+n}k b^k z_1^k d^{l+n-k} z_2^{l+n-k}.
\end{align*}
Hence the \LHS\ of \eqref{2} can be rewritten as
\begin{equation}
\binom{2l}{l-n}^\half\, \sum_{j=0}^{l-n} \sum_{k=0}^{l+n}
\binom{l-n}j \binom{l+n}k
a^j b^k c^{l-n-j} d^{l+n-k} z_1^{j+k} z_2^{2l-j-k}.
\label{14}
\end{equation}
In this double sum we make a change of summation variables
$(j,k)\mapsto(m,j)$, where
$j+k=l-m$. Hence
\begin{equation}
(j,k) \mapsto (l-k-j,j)\;\;
\mbox{with inverse map}\;\;
(m,j)\mapsto(j,l-m-j).
\label{12}
\end{equation}
Now we have
\begin{multline}
0\le j\le l-n\;\;{\rm and}\;\;0\le k\le l+n\;\;
\Equiva\\
-l\le m\le l\;\;{\rm and}\;\;0\le j\le l-n\;\;{\rm and}\;\;
-m-n\le j\le l-m.
\label{11}
\end{multline}
Indeed, the inequalities to the left of the equivalence sign in \eqref{11}
imply that $0\le j+k\le 2l$, hence $0\le l-m\le 2l$, hence $-l\le m\le l$.
Also, $0\le k\le l+n$ implies $0\le l-m-j\le l+n$, hence
$-m-n\le j\le l-m$.
Conversely, $-m-n\le j\le l-m$ implies (substitute $m=l-k-j$) that
$-l-n+k+j\le j\le k+j$, hence $0\le k\le l+n$.
(Note that $-l\le m\le l$ to the right of the equivalence sign in \eqref{11}
is not strictly needed
because it is implied by the other inequaltities on the right.)

We conclude that the double sum \eqref{14} can be rewritten by the
substitution $j+k=l-m$ as follows:
\begin{equation}
\binom{2l}{l-n}^\half \sum_{m=-l}^l\; \sum_{j=0\vee(-m-n)}^{(l-m)\wedge(l-n)}
\binom{l-n}j \binom{l+n}{l-m-j}
a^j b^{l-m-j} c^{l-n-j} d^{n+m+j} z_1^{l-m} z_2^{l+m}.
\label{13}
\end{equation}
Here the first summation is by convention over all $m\in\{-l,-l+1,\ldots,l\}$.
In the second summation the symbol $\vee$ means maximum and the symbol
$\wedge$ means minimum. The range of the double summation in \eqref{13}
is justified by the equivalence \eqref{11}. Note that the second summation
is an inner summation since its summation bounds depend on $m$, which is the
summation variable for the outer summation.
The summand in \eqref{13} is obtained from the summand in 
\eqref{14} by the substitution $k=l-m-j$.

Since \eqref{13} is a rewritten form of the \LHS\ of \eqref{2}, it must
be equal to the \RHS\ of \eqref{2}. Both \eqref{13} and the
\RHS\ of \eqref{2} are polynomials in $z_1,z_2$ with explicit coefficients.
Hence the corresponding coefficients must be equal. We conclude:
\begin{proposition}
\begin{equation}
t_{m,n}^l\bma a&b\\c&d\ema=
\binom{2l}{l-m}^{-\half} \binom{2l}{l-n}^\half
\sum_{j=0\vee(-m-n)}^{(l-m)\wedge(l-n)}
\binom{l-n}j \binom{l+n}{l-m-j}
a^j b^{l-m-j} c^{l-n-j} d^{n+m+j}.
\label{15}
\end{equation}
\end{proposition}

Note that the summation bounds in (1.12) reduce to one of four alternatives
depending on the signs of $m+n$ and $m-n$:
\begin{align*}
0\le j\le l-m\quad&{\rm if}\quad m+n\ge0\;\;{\rm and}\;\;m-n\ge0;\\
0\le j\le l-n\quad&{\rm if}\quad m+n\ge0\;\;{\rm and}\;\;m-n\le0;\\
-m-n\le j\le l-m\quad&{\rm if}\quad m+n\le0\;\;{\rm and}\;\;m-n\ge0;\\
-m-n\le j\le l-n\quad&{\rm if}\quad m+n\le0\;\;{\rm and}\;\;m-n\le0.
\end{align*}
These four alternatives correspond two four subsets of the set
$\{(m,n)\mid m,n\in\{-l,-l+1,\ldots,l\}\}$, which have triangular shape,
overlapping boundaries, and together span the whole set. These four subsets
are mapped onto each other by the symmetries \eqref{5}--\eqref{7}.

Hence it is sufficient to compute $t_{m,n}^l$ if $m+n\ge0$, $m-n\ge0$.
For a while we only assume $m+n\ge0$ and not yet $m-n\ge0$
Then \eqref{15} takes the form
\begin{equation}
t_{m,n}^l\bma a&b\\c&d\ema=
\binom{2l}{l-m}^{-\half} \binom{2l}{l-n}^\half\;
\sum_{j\ge0}
\binom{l-n}j \binom{l+n}{l-m-j}
a^j b^{l-m-j} c^{l-n-j} d^{n+m+j}.
\label{16}
\end{equation}

We will rewrite the \RHS\ of \eqref{16} first as a Gauss hypergeometric
function (with some elementary factors in front)
and next as a Jacobi polynomial.
For this derivation remember the {\em Pochhammer symbol}
\begin{equation*}
(a)_0:=1,\qquad
(a)_k:=a(a+1)\ldots(a+k-1)\quad(k=1,2,\ldots).
\end{equation*}
In particular, note that
\begin{equation*}
\frac{(n+k)!}{n!}=(n+1)_k\,,\qquad
\frac{n!}{(n-k)!}=(-1)^k (-n)_k\,.
\end{equation*}
Now we have
\begin{align}
&\binom{2l}{l-m}^{-\half} \binom{2l}{l-n}^\half\;
\sum_{j\ge0}
\binom{l-n}j \binom{l+n}{l-m-j}
a^j b^{l-m-j} c^{l-n-j} d^{n+m+j}
\nonu\\
&=\left(\frac{(l+m)!\,(l-m)!}{(l+n)!\,(l-n)!}\right)^\half\,
\sum_{j\ge0}\frac{(l-n)!}{j!\,(l-n-j)!}\,
\frac{(l+n)!}{(l-m-j)!\,(n+m+j)!}\,
a^j b^{l-m-j} c^{l-n-j} d^{n+m+j}
\nonu\\
&=\left(\frac{(l+m)!\,(l-m)!}{(l+n)!\,(l-n)!}\right)^\half\,
\frac{(l+n)!\,b^{l-m} c^{l-n} d^{m+n}}{(l-m)!\,(m+n)!}\,
\sum_{j\ge0} \frac{(l-m)!}{(l-m-j)!}\,\frac{(l-n)!}{(l-n-j)!}\,
\frac{(m+n)!}{(m+n+j)!\,j!}\,\left(\frac{ad}{bc}\right)^j
\nonu\\
&=\left(\frac{(l+m)!\,(l+n)!}{(l-m)!\,(l-n)!}\right)^\half\,
\frac{b^{l-m} c^{l-n} d^{m+n}}{(m+n)!}\,
\sum_{j\ge0}\frac{(-l+m)_j\,(-l+n)_j}{(m+n+1)_j\,j!}\,
\left(\frac{ad}{bc}\right)^j
\nonu\\
&=\left(\frac{(l+m)!\,(l+n)!}{(l-m)!\,(l-n)!}\right)^\half\,
\frac{b^{l-m} c^{l-n} d^{m+n}}{(m+n)!}\,
\hyp21{-l+m,-l+n}{m+n+1}{\frac{ad}{bc}}.
\label{17}
\end{align}
Here we used the definition of Gauss hypergeometric series, see
\cite[Chapter 2]{2}.
Note that the two upper parameters $-l+m,-l+n$ of the hypergeometric function
in \eqref{17} are both non-positive, and that the series will terminate after
the term with $j=(l-m)\wedge(l-n)$.

Pfaff's transformation (see \cite[(2.2.6)]{2})
\begin{equation}
\hyp21{a,b}cz=(1-z)^{-a}\,\hyp21{a,c-b}c{\frac z{z-1}}
\label{44}
\end{equation}
implies for the hypergeometric function
in \eqref{17} that
\begin{equation*}
\hyp21{-l+m,-l+n}{m+n+1}{\frac{ad}{bc}}=
b^{m-l}c^{m-l}(bc-ad)^{l-m}\hyp21{-l+m,l+m+1}{m+n+1}{\frac{ad}{ad-bc}}.
\end{equation*}
Hence we arrive at the following rewritten form of \eqref{16}
(from now on assume $m+n\ge0$, $m-n\ge0$):
\begin{equation}
t_{m,n}^l\bma a&b\\c&d\ema=
\left(\frac{(l+m)!\,(l+n)!}{(l-m)!\,(l-n)!}\right)^\half\,
\frac{c^{m-n} d^{m+n} (bc-ad)^{l-m}}{(m+n)!}\,
\hyp21{-l+m,l+m+1}{m+n+1}{\frac{ad}{ad-bc}}.
\label{18}
\end{equation}

Now use the expression of {\em Jacobi polynomials} in terms of the
Gauss hypergeometric function (see \cite[Definition 2.5.1]{2}):
\begin{align}
P_n^{(\al,\be)}(x)&=\frac{(\al+1)_n}{n!}\,
\hyp21{-n,n+\al+\be+1}{\al+1}{\frac{1-x}2}
\label{38}\\
&=\sum_{k=0}^n
\frac{(n+\al+\be+1)_k (\al+k+1)_{n-k}}{k!\,(n-k)!}\,
\Big(\frac{x-1}2\Big)^k.
\label{48}
\end{align}
Although Jacobi polynomials are usually considered for $\al,\be>-1$
because of a nice orthogonality property \eqref{49}
for these parameter values, they depend polynomially on $\al,\be$
by \eqref{48}, and are therefore well-defined for all $\al,\be$.

Formula \eqref{38}
implies for the hypergeometric function in \eqref{18} that
\begin{equation*}
\hyp21{-l+m,l+m+1}{m+n+1}{\frac{ad}{ad-bc}}
=\frac{(l-m)!\,(m+n)!}{(l+n)!}\,
P_{l-m}^{(m+n,m-n)}\left(\frac{bc+ad}{bc-ad}\right).
\end{equation*}
Hence we can further rewrite \eqref{18} (if $m+n\ge0$, $m-n\ge0$) as follows:
\begin{equation}
t_{m,n}^l\bma a&b\\c&d\ema=
\left(\frac{(l+m)!\,(l-m)!}{(l+n)!\,(l-n)!}\right)^\half\,
c^{m-n} d^{m+n} (bc-ad)^{l-m}\,
P_{l-m}^{(m+n,m-n)}\left(\frac{bc+ad}{bc-ad}\right).
\label{19}
\end{equation}
We are in particular interested in \eqref{19} if
$\bma a&b\\c&d\ema\in SU(2)$. Note that by \eqref{20} a general element
of SU(2) can be written as
\begin{equation*}
\bma \sin\tha\,e^{i\phi}&-\cos\tha\,e^{-i\psi}\\
\cos\tha\,e^{i\psi}&\sin\tha\,e^{-i\phi} \ema\qquad
\mbox{with $0\le\tha\le\pi/2$ and $\phi,\psi\in[0,2\pi)$.}
\end{equation*}
Hence we obtain:
\begin{theorem}
If $m+n\ge0$, $m-n\ge0$ then
\begin{multline}
t_{m,n}^l\bma \sin\tha\,e^{i\phi}&-\cos\tha\,e^{-i\psi}\\
\cos\tha\,e^{i\psi}&\sin\tha\,e^{-i\phi} \ema=
(-1)^{l-m}
\left(\frac{(l+m)!\,(l-m)!}{(l+n)!\,(l-n)!}\right)^\half\\
\times
e^{-i(m+n)\phi}\,e^{i(m-n)\psi}\,
(\sin\tha)^{m+n} (\cos\tha)^{m-n}\,
P_{l-m}^{(m+n,m-n)}(\cos 2\tha).
\label{21}
\end{multline}
\end{theorem}
\section{A Rodrigues type formula for the matrix elements}
From the generating function \eqref{2} we can derive a quite
different explicit formula for the matrix elements $t_{m,n}^l$.
In \eqref{2} put
\[
\bma a&b\\c&d\ema:=\bma \sin\tha&-\cos\tha\\ \cos\tha&\sin\tha\ema,
\quad
z_1:=\thalf(s-\cos 2\tha),\quad
z_2:=\sin\tha\cos\tha.
\]
Then \eqref{2} takes the form
\begin{multline*}
\binom{2l}{l-n}^\half (\sin\tha)^{l-n} (\cos\tha)^{l+n}\,
(1-s)^{l-n} (1+s)^{l+n}\\
=\sum_{m=-l}^l \binom{2l}{l-m}^\half (2\sin\tha\cos\tha)^{l+m}
t_{m,n}^l\bma \sin\tha&-\cos\tha\\ \cos\tha&\sin\tha\ema
(s-\cos 2\tha)^{l-m}.
\end{multline*}
Hence
\begin{multline}
t_{m,n}^l\bma \sin\tha&-\cos\tha\\ \cos\tha&\sin\tha\ema
=\Bigg(\frac{(l+m)!}{(l-m)!\,(l+n)!\,(l-n)!}\Bigg)^\half
2^{-l-m}
(\sin\tha)^{-m-n} (\cos\tha)^{-m+n}\\
\times \Big(\frac d{ds}\Big)^{l-m}\Big((1-s)^{l+n}(1+s)^{l-n}\Big)
\Big|_{s=\cos2\tha}\,.
\label{50}
\end{multline}
In view of the Rodrigues formula (see \cite[(2.5.13\!\'{})]{2})
\begin{equation}
P_n^{(\al,\be)}(x)=\frac{(-1)^n}{2^n n!}\,(1-x)^{-\al}(1+x)^{-\be}
\Big(\frac d{dx}\Big)^n \Big((1-x)^{n+\al}(1+x)^{n+\be}\Big),
\label{51}
\end{equation}
valid for arbitrary $\al,\be$, formula \eqref{50} implies
formula \eqref{21} for $\phi=\psi=0$, not just for $m+n,m-n\ge0$,
but for all $m,n\in\{-l,-l+1,\ldots,l\}$.
Conversely, \eqref{21} and \eqref{50} together imply \eqref{51}
for $\al,\be\in\{0,1,2,\ldots\}$.
\section{Orthogonality of matrix elements}
We introduce a special Borel measure $\mu$ on SU(2) such that
\begin{equation}
\int_{\rm SU(2)} f d\mu=
\frac1{2\pi^2}\int_{\phi=0}^{2\pi} \int_{\psi=0}^{2\pi} \int_{\tha=0}^{\pi/2}
f\bma \sin\tha\,e^{i\phi}&-\cos\tha\,e^{-i\psi}\\
\cos\tha\,e^{i\psi}&\sin\tha\,e^{-i\phi} \ema\,\sin\tha\,\cos\tha\
d\tha\,d\psi\,d\phi.
\label{22}
\end{equation}
for all continuous funtions $f$ on SU(2).
Note that
\begin{equation}
\int_{\rm SU(2)} d\mu=1.
\label{23}
\end{equation}

The matrix elements $t_{m,n}^l$ satisfy a remarkable orthogonality
relation with respect to this measure:
\begin{equation}
\int_{SU(2)} t_{m,n}^l\,\overline{t_{m',n'}^{l'}}\,d\mu=
\frac1{2l+1}\,\de_{l,l'}\,\de_{m,m'}\,\de_{n,n'}\,.
\label{24}
\end{equation}
For $(m,n)\ne(m',n')$ this follows immediately from \eqref{21}, \eqref{22}
and the symmetries \eqref{5}--\eqref{7}.
(exercise). For $(m,n)=(m',n')$ with $m+n,m-n\ge0$
we have to show that
\begin{multline*}
\frac{(l+m)!\,(l-m)!}{(l+n)!\,(l-n)!}\,
\int_0^{\pi/2} P_{l-m}^{(m+n,m-n)}(\cos 2\tha)\,
P_{l'-m}^{(m+n,m-n)}(\cos 2\tha)\\
\times(\sin\tha)^{2m+2n+1} (\cos\tha)^{2m-2n+1}\,d\tha
=\frac1{2l+1}\,\de_{l,l'}\,.
\end{multline*}
By the substitution $x=\cos2\tha$ this can be rewritten as
\begin{equation}
\frac{(l+m)!\,(l-m)!}{(l+n)!\,(l-n)!\,2^{2m+1}}\,
\int_{-1}^1 P_{l-m}^{(m+n,m-n)}(x)\,
P_{l'-m}^{(m+n,m-n)}(x)\,(1-x)^{m+n} (1+x)^{m-n}\,dx=
\frac{\de_{l,l'}}{2l+1}\,.
\label{25}
\end{equation}
In order to show this identity we remember the orthogonality relations
for Jacobi polynomials (see \cite[(2.5.14)]{2}):
\begin{equation}
\int_{-1}^1 P_m^{(\al,\be)}(x)\,P_n^{(\al,\be)}(x)\,(1-x)^\al (1+x)^\be\,dx
=h_n^{(\al,\be)}\,\de_{m,n}
\label{49}
\end{equation}
with
\begin{equation}
h_n^{(\al,\be)}=
\frac{2^{\al+\be+1}(n+\al+\be+1)_n\,\Ga(n+\al+1)\,\Ga(n+\be+1)}
{n!\,\Ga(2n+\al+\be+2)}\,.
\end{equation}
Now observe that
\begin{equation*}
h_{l-m}^{(m+n,m-n)}=
\frac{2^{2m+1} (l+n)!\,(l-n)!}{(2l+1)(l+m)!\,(l-m)!}\,.
\end{equation*}
This settles \eqref{25} and hence \eqref{24}.
\bPP
The orthogonality relation \eqref{24} is a special case of
{\em Schur's orthogonality relations} for the matrix elements of the
irreducuible unitary representations of a compact group. For the
formulation of this theorem we need the concept of the {\em Haar measure}
(see for instance \cite[Section 15]{17})
on a compact group.
\begin{theorem}
Let $G$ be a compact group. There is a unique Borel measure $\mu$ on $G$,
called {\em Haar measure}, such that $\mu(G)=1$ and,
for all Borel sets $E\subset G$ and for all $g\in G$,
\begin{equation*}
\mu(g E) =\mu(E) =\mu(E g).
\end{equation*}
\end{theorem}

For a continuous function $f$ on $G$ this implies:
\begin{equation*}
\int_G f(hg)\,d\mu(g)=\int_G f(g)\,d\mu(g)=\int_G f(gh)\,d\mu(g)\qquad(h\in G).
\end{equation*}

We also need the concept of equivalence of representations:
\begin{definition}
Let $\pi_1$ and $\pi_2$ be representations of a group $G$ on finite
dimensional complex vector spaces $V_1$ and $V_2$, respectively.
Then $\pi_1$ and $\pi_2$ are called {\em equivalent} to each other
if there is a linear bijection $A\colon V\to W$ such that
\begin{equation*}
A\,\pi_1(g)=\pi_2(g)\,A\quad\mbox{for all $g\in G$.}
\end{equation*}
The representations $\pi_1$ and $\pi_2$ are called {\em inequivalent} to
each other if they are not equivalent to each other.
\end{definition}

The relation of equivalence of representations is an equivalence relation
on the set of all finite dimensional representations of a group $G$
(see for instance \cite [Theorem 27.19]{18}).
\begin{theorem}
\label{26}
{\bf(Schur's orthogonality relations)}\\
Let $G$ be a compact group with Haar measure $\mu$. Let $\pi$ and $\rho$
be finite dimensional complex
irreducible unitary representations of $G$ which are inequivalent
to each other. Let $\pi$ resp.\ $\rho$ have matrix elements
$(\pi_{i,j})_{i,j=1,\ldots,d_\pi}$ and $(\rho_{k,l})_{k,l=1,\ldots,d_\rho}$
with respect to certain orthonormal bases of their representation
spaces. Then
\begin{equation*}
\int_G \pi_{i,j}(g)\,\overline{\rho_{k,l}(g)}\,d\mu(g)=0
\end{equation*}
and
\begin{equation*}
\int_G \pi_{i,j}(g)\,\overline{\pi_{k,l}(g)}\,d\mu(g)
=\frac 1{d_\pi}\,\de_{i,k}\,\de_{j,l}.
\end{equation*}
\end{theorem}

Let $(\pi^{\al})_{\al\in A}$ be a maximal set of mutually inequivalent
finite dimensional complex irreducible unitary representations of $G$,
and put $d_\al:=d_{\pi^\al}$.
Then the orthogonality relations in Theorem \ref{26} can be written as
\begin{equation}
\int_G \pi_{i,j}^\al(g)\,\overline{\pi_{k,l}^\be(g)}\,d\mu(g)=
\frac1{d_\al}\,\de_{\al,\be}\,\de_{i,k}\,\de_{j,l}.
\label{27}
\end{equation}
Write $L^2(G)$ for $L^2(G,\mu)$. The functions $\pi_{i,j}$ are continuous
on $G$, so they are certainly in $L^2(G)$. By \eqref{27} the functions
$d_\al^\half\,\pi_{i,j}^\al$ ($\al\in A$, $i,j=1,\ldots,d_\al$)
form an orthonormal system in $L^2(G)$. Then the {\em Peter-Weyl theorem}
(see for instance \cite[Theorem 27.49]{18}) says:
\begin{theorem}
With notation as above, the functions $d_\al^\half\,\pi_{i,j}^\al$
form a complete orthonormal system in $L^2(G)$.
\end{theorem}

In order to match the orthogonality relation \eqref{24} to \eqref{27} we
need an explicit form of the Haar measure on SU(2) and we have to prove
that the representations $t^l$ are irreducible.
First we deal with the Haar measure.
By \eqref{20} the group SU(2) is homeomorphic with the unit sphere
$S^3=\{(a,c)\in \CC^2\mid |a|^2+|c|^2=1\}$.
Let $S\in {\rm SU(2)}$.
A left multiplication $T\mapsto ST\colon {\rm SU(2)}\to{\rm SU(2)}$
corrresponds to some
rotation of $S^3$ (exercise). Thus a rotation invariant measure on $S^3$
will provide, after suitable normalization, the Haar measure on SU(2).
There exists, up to a constant factor,
a unique rotation invariant measure $\om$ on $S^3$. This measure is
such that, for all continuous functions on $\RR^4$ of compact support
and with $\la$ Lebesgue measure on $\RR^4$,
\begin{equation}
\int_{\RR^4} f\,d\la=
\int_{r=0}^\iy \int_{\xi\in S^3} f(r\xi)\,d\om(\xi)\,r^3\,dr.
\label{28}
\end{equation}
Now take coordinates
\begin{equation*}
x=(r\sin\tha\cos\phi,r\sin\tha\sin\phi,r\cos\tha\cos\psi,r\cos\tha\sin\psi)
\end{equation*}
on $\RR^4$, which means for $r=1$ that
$x_1+ix_2=\sin\tha\,e^{i\phi}$, $x_3+ix_4=\cos\tha\,e^{i\psi}$.
These are just the coordinates chosen in \eqref{22} for $(a,c)\in\CC^2$
with $|a|^2+|c|^2=1$. A straightforward computation of the Jacobian yields:
\begin{multline}
\int_{\RR^4} f(x_1,x_2,x_3,x_4)\,dx_1\,dx_2\,dx_3\,dx_4\\
=
\int_{r=0}^\iy \int_{\tha=0}^{\pi/2} \int_{\phi=0}^{2\pi} \int_{\psi=0}^{2\pi}
f(r\sin\tha\cos\phi,r\sin\tha\sin\phi,r\cos\tha\cos\psi,r\cos\tha\sin\psi)\\
\times r^3 \sin\tha \cos\tha\,dr\,d\tha\,d\phi\,d\psi.
\label{29}
\end{multline}
Comparison of \eqref{28} and \eqref{29} gives, for continuous functions $F$
on $S^3\subset\CC^2$, that
\begin{equation}
\int_{S^3} F\,d\om=
\int_{\tha=0}^{\pi/2} \int_{\phi=0}^{2\pi} \int_{\psi=0}^{2\pi}
F(\sin\tha\,e^{i\phi},\cos\tha\,e^{i\psi})\,
\sin\tha \cos\tha\,d\tha\,d\phi\,d\psi.
\end{equation}
In view of the previous observations we have thus shown that the Haar measure
on SU(2) is given by \eqref{22}.
\mPP
Now we will show the irreducibility of the representations $t^l$ as
representations of SU(2).

Put $a_\phi:=\bma e^{i\phi}&0\\0&e^{-i\phi}\ema$. Then
$a_\phi a_\psi=a_{\phi+\psi}$ and $a_{\phi+2\pi}=a_\phi$.
The group $A:=\{a_\phi\mid 0\le\phi<2\pi\}$ is a closed abelian subgroup
of SU(2). It is isomorphic and homeomorphic with the group $U(1)$
of complex numbers of absolute value 1, which has multiplication of
complex numbers as the group multiplication.
It follows from \eqref{31} that
\begin{equation}
t^l(a_\phi)\,\psi_n^l=e^{-2in\phi}\,\psi_n^l.
\label{33}
\end{equation}
\begin{lemma}
Let $V$ be an invariant subspace of $\FSH_l$ with respect to the representation
$\pi^l$ of SU(2).
If $v\in V$ and $\lan v,\psi_m^l\ran\ne0$ then $\psi_m^l\in V$.
\end{lemma}
\Proof
We have
\begin{align*}
v&=\sum_{n=-l}^l \lan v,\psi_n^l\ran\,\psi_n^l,\\
t^l(a_\phi)\,v&=\sum_{n=-l}^l \lan v,\psi_n^l\ran\,t^l(a_\phi)\,\psi_n^l
=\sum_{n=-l}^l \lan v,\psi_n^l\ran\,e^{-2in\phi}\,\psi_n^l.\\
\noalign{\mbox{Hence}}
\int_0^{2\pi}e^{2im\phi}\,t^l(a_\phi)\,v\,d\phi&=
\sum_{n=-l}^l \lan v,\psi_n^l\ran
\left(\int_0^{2\pi}e^{2im\phi} e^{-2in\phi}\,d\phi\right)\psi_n^l
=2\pi \lan v,\psi_m^l\ran \psi_m^l.
\end{align*}
The integral on the left should be interpreted as a Riemann integral of
vectors, which can be approximated by Riemann sums of vectors.
Since $v\in V$, each approximating Riemann sum is in $V$, and hence also
their limit, the Riemann integral, is in $V$.
Hence $2\pi \lan v,\psi_m^l\ran \psi_m^l\in V$. So $\psi_m^l\in V$ if
$\lan v,\psi_m^l\ran\ne0$.
\qed
\bPP
This Lemma implies the following Proposition, the proof of which is left
as an exercise.
\begin{proposition}
Let $V$ be an invariant subspace of $\FSH_l$
with respect to the representation $\pi^l$ of {\rm SU(2)}.
Then there is a subset
$\FSA$ of $\{-l,\ldots,l\}$ such that $V=\Span\{\psi_n^l\mid n\in\FSA\}$.
Let $W$ be the orthoplement of $V$ and $\FSB$ the complement of $\FSA$.
Then $W$ is also an invariant subspace and
$W=\Span\{\psi_n^l\mid n\in\FSB\}$.
\label{34}
\end{proposition}
\begin{theorem}
The representation $t^l$ of {\rm SU(2)} is irreducible.
\end{theorem}
\Proof
Suppose $t^l$ is not irreducible.
By Proposition \ref{34} $\FSH_l$ is the orthogonal direct sum of
invariant subspaces
$V=\Span\{\psi_n^l\mid n\in\FSA\}$ and $W=\Span\{\psi_n^l\mid n\in\FSB\}$,
where $\{-l,\ldots,l\}$ is the disjoint union of certain nonempty subsets
$\FSA$ and $\FSB$. One of these subsets, say $\FSA$, will contain $l$.
Then some $m$ will be in $\FSB$. Then $t^l(T)\,\psi_l^l$ will be in $V$
for all $T\in {\rm SU(2)}$, and therefore orthogonal to $\psi_m^l$.
Hence $t_{m,l}^l(T)=0$ for all $T\in {\rm SU(2)}$.
In particular, also using \eqref{35}, we obtain
\begin{equation*}
0=t_{m,l}^l\bma\sin\tha&-\cos\tha\\\cos\tha&\sin\tha\ema=
(-1)^{l-m}\,(\cos\tha)^{l-m}\,(\sin\tha)^{l+m},
\end{equation*}
which gives a contradiction.
\qed
\bLP
{\bf Remark}\quad
The above proof of irredicibility is by the so-called global
(or non-infinitesimal) method, cf.~\cite{20}.
See for instance \cite[Section III.2.3]{12} for a proof using the
infinitesimal method, i.e., considering the corresponding representation
of the Lie algebra of SU(2).
\bPP
So finally we have matched \eqref{24} to \eqref{27}.
We could have started with \eqref{27} and have derived from this \eqref{24}
and hence \eqref{25}. Thus a proof of the orthogonality relations for
Jacobi polynomials with nonnegative integer parameters is possible
from the interpretation of Jacobi polynomials in connection with SU(2).
\section{Exercises}
\label{53}
{\bf 1.}
(An interpretation of Krawtchouk polynomials as matrix elements of
irreducible representations of SU(2), see \cite[Section 2]{19})
\sLP
{\bf a)}\quad
Prove that
\begin{equation}
P_n^{(\al,\be)}(-x)=(-1)^n\,P_n^{(\be,\al)}(x).
\label{36}
\end{equation}
(Use the orthogonality relations for Jacobi polynomials and the explicit
expression for the leading coefficient of $P_n^{(\al,\be)}(x)$.)
\sLP
{\bf b)}\quad
Prove that
\begin{equation}
\hyp21{-n,b}cx=\frac{(c-b)_n}{(c)_n}\,
\hyp21{-n,b}{b-c-n+1}{1-x}\qquad(n=0,1,2,\ldots).
\label{37}
\end{equation}
(Use \eqref{38} and \eqref{36}.)
\sLP
{\bf c)}\quad
Prove that
\begin{equation}
\hyp21{-n,-m}cx=\frac{(c)_{m+n}}{(c)_n(c)_m}\,
\hyp21{-n,-m}{-c-n-m+1}{1-x}\qquad(n,m=0,1,2,\ldots).
\label{39}
\end{equation}
(Use \eqref{37}.)
\sLP
{\bf d)}\quad
Prove that, for $m+n\ge0$,
\begin{equation}
t_{m,n}^l\bma a&b\\c&d\ema=
\binom{2l}{l-m}^\half
\binom{2l}{l-n}^\half
b^{l-m} c^{l-n} d^{m+n}\,
\hyp21{-l+m,-l+n}{-2l}{\frac{bc-ad}{bc}}.
\label{40}
\end{equation}
(Use \eqref{39} and \eqref{17}.)
\sLP
{\bf e)}\quad
Prove that, for $m+n\ge0$,
\begin{multline}
t_{m,n}^l\bma\sin\tha&-\cos\tha\\\cos\tha&\sin\tha\ema=
\binom{2l}{l-m}^\half
\binom{2l}{l-n}^\half
(-1)^{l-m} (\cos\tha)^{2l-m-n} (\sin\tha)^{m+n}\\
\times K_{l-m}(l-n;\cos^2\tha,2l),
\label{41}
\end{multline}
where the {\em Krawtchouk polynomials} are given by
\begin{equation}
K_n(x;p,N):=\hyp21{-n,-x}{-N}{p^{-1}}\qquad
(n=0,1,\ldots,N).
\label{42}
\end{equation}
\sLP
{\bf f)}\quad
Prove that
\begin{equation}
K_n(x;p,N)=(1-p^{-1})^{x+n-N}\,K_{N-n}(N-x;p,N).
\label{43}
\end{equation}
(Use \eqref{42} and Euler's transformation for Gauss hypergeometric functions.)
\sLP
{\bf g)}\quad
Prove that \eqref{41} remains valid for all $m,n$.\\
(Use \eqref{6} and \eqref{43}.)
\sLP
{\bf h)}\quad
Show that
\begin{equation}
\sum_{n=-l}^l t_{m,n}^l\bma\sin\tha&-\cos\tha\\\cos\tha&\sin\tha\ema
t_{m',n}^l\bma\sin\tha&-\cos\tha\\\cos\tha&\sin\tha\ema=\de_{m,m'}\,,
\label{45}
\end{equation}
and that this matches with the orthogonality relation for the
Krawtchouk polynomials occurring on the \RHS\ of \eqref{41}.
\bLP
{\bf 2.} (Addition formula and product formula
for Legendre polynomials)\quad
Let $l=0,1,2,\ldots\;$.
\sLP
{\bf a)}\quad
Prove that, for $ad-bc=1$,
\begin{equation}
t_{0,0}^l\bma a&b\\c&d\ema=P_l(2ad-1),
\label{46}
\end{equation}
where $P_l$ is the Legendre polynomial.
(Use \eqref{19}.)
\sLP
{\bf b)}\quad
Prove that
\begin{equation}
t_{0,0}^l(T)=P_l(\cos\tha_1\cos\tha_2+\sin\tha_1\sin\tha_2\cos\phi)
\end{equation}
if
\begin{equation}
T=\bma\sin\thalf\tha_1&-\cos\thalf\tha_1\sLP
\cos\thalf\tha_1&\sin\thalf\tha_1\ema
\bma e^{\half i\phi}&0\\0&e^{-\half i\phi}\ema
\bma\sin\thalf\tha_2&\cos\thalf\tha_2\sLP
-\cos\thalf\tha_2&\sin\thalf\tha_2\ema.
\end{equation}
{\bf c)}\quad
Prove that
\begin{multline}
P_l(\cos\tha_1\cos\tha_2+\sin\tha_1\sin\tha_2\cos\phi)=
P_l(\cos\tha_1)\,P_l(\cos\tha_2)\\
+\sum_{0<|k|\le l} t_{0,k}^l
\bma\sin\thalf\tha_1&-\cos\thalf\tha_1\sLP
\cos\thalf\tha_1&\sin\thalf\tha_1\ema
t_{k,0}^l
\bma\sin\thalf\tha_2&\cos\thalf\tha_2\sLP
-\cos\thalf\tha_2&\sin\thalf\tha_2\ema
e^{-ik\phi}.
\end{multline}
{\bf d)}\quad
Prove that
\begin{equation}
P_l(\cos\tha_1)\,P_l(\cos\tha_2)=
\frac1{2\pi}
\int_0^{2\pi}
P_l(\cos\tha_1\cos\tha_2+\sin\tha_1\sin\tha_2\cos\phi)\,d\phi.
\end{equation}
\section{Some history}
\label{52}
The irreducible matrix representation $t^l$ of SU(2),
introduced in Section \ref{47},
is known as the {\em Wigner d-matrix}. Wigner introduced this in
\cite[(10)]{1} (1927), but only for representations of ${\rm SO}(3)$, i.e.,
for $l=0,1,2,\ldots\;$. In a series of papers \cite{15} (1928)
together with von Neumann half integer values of $l$, i.e., double-valued
representations of SO(3) were also allowed
in order to accommodate electron spin.
Next, in Chapter 15 of his book
\cite{5} (1931; see also a translated and extended version
\cite{6} (1959)),
Wigner observes the two-to-one homomorphism ${\rm SU(2)} \to {\rm SO(3)}$
and he obtains the explicit formula \eqref{15}
for the matrix elements of the irreducible representations
$t^l$ ($l=0,\thalf,1,\tfrac32,\ldots$) of SU(2),
but he does not express this in terms of a
hypergeometric function or a Jacobi polynomial.
This is also the case in Talman's book \cite[Section 8,2]{14},
which is based on lectures by Wigner.

G\"uttinger, in the same year 1931 as \cite{5},
in a {\em Mathematischer Anhang}
to his paper~\cite{10} (Diplomarbeit at TH Z\"urich under Pauli),
expresses the matrix elements of the irreducible representations
of SU(2) in terms of Jacobi polynomials, on the one hand via the
Rodrigues formula (see \eqref{50}, \eqref{51}),
on the other hand
by recognizing the matrix elements as terminating Gauss hypergeometric
series (see \eqref{17}).
Quite remarkably, G\"uttinger's paper does not have any reference
to work by Wigner.
In 1952 Schwinger \cite[Section 2]{11} and
Gelfand \& \v{S}apiro \cite[Section 7.4]{7}
(also part I of the book \cite{9}) give similar treatments, but only
Schwinger refers to \cite{10}, 
while \cite{7} has no references at all (but in the bibliography in the book
\cite{9} Wigner's books \cite{5}, \cite{6} are present).
Gelfand \& Graev \cite[Section 1]{16} (1965) briefly repeat the material
in \cite{7} about representations of SU(2) and Jacobi polynomials
before they pass to representations of GL(n) for $n>2$.
Anyhow, the connection of Jacobi polynomials with representation theory
of SU(2) seemed to be common knowledge in 1965, since
Coleman, in his review MR0201568 in Math.\ Reviews
of this paper \cite{16},
writes: ``It is well-known that explicit expressions are available for
the finite-dimensional irreducible representations of the unitary group
in two dimensions, SU(2), and of the full linear group, GL(2,C).
The representation spaces may be realized conveniently as spaces of polynomials in a single variable and their matrix elements involve
Jacobi polynomials.'' In the same year 1965 the original Russian edition
of Vilenkin's book \cite{12} appears, with extensive treatment of this
theory in Chapter~III.
\quad\\
\begin{footnotesize}
\begin{quote}
{ T. H. Koornwinder, Korteweg-de Vries Institute, University of
 Amsterdam,\\
 P.O.\ Box 94248, 1090 GE Amsterdam, The Netherlands;

\vspace{\smallskipamount}
email: }{\tt T.H.Koornwinder@uva.nl}
\end{quote}
\end{footnotesize}

\end{document}